\documentclass[12pt]{amsart}

\usepackage{amssymb, amsmath, amscd}

\usepackage{mathptmx}

\newtheorem{teo}{\bf Theorem}[section]
\newtheorem{thm}[teo]{\bf Theorem}
\newtheorem{lem}[teo]{\bf Lemma}
\newtheorem{cor}[teo]{\bf Corollary}
\newtheorem{prop}[teo]{\bf Proposition}

\newtheorem{quest}[teo]{\bf Question}

\theoremstyle{definition}
\newtheorem{defi}[teo]{\bf Definition}

\newtheorem{rem}[teo]{\bf Remark}
\newtheorem{exa}[teo]{\bf Example}

\newcommand{\Z}{\mathbb{Z}}
\newcommand{\R}{\mathbb{R}}
\newcommand{\N}{\mathbb{N}}

\DeclareMathOperator{\GL}{\mathrm{GL}}
\DeclareMathOperator{\Tor}{\mathrm{Tor}}
\DeclareMathOperator{\reg}{\mathrm{reg}}
\DeclareMathOperator{\revlex}{\mathrm{revlex}}

\DeclareMathOperator{\p}{\mathfrak{p}}
\DeclareMathOperator{\m}{\mathfrak{m}}

\DeclareMathOperator{\eb}{\bf e}

\DeclareMathOperator{\vb}{\bf v}

\DeclareMathOperator{\gin}{gin}
\DeclareMathOperator{\gins}{gins}
\DeclareMathOperator{\lex}{lex}

\DeclareMathOperator{\shift}{\Delta}
\DeclareMathOperator{\facets}{facets}
\DeclareMathOperator{\init}{init}
\DeclareMathOperator{\adeg}{adeg}

\DeclareMathOperator{\Soc}{Soc}

\DeclareMathOperator{\ini}{in}

\DeclareMathOperator{\Deg}{Deg} \DeclareMathOperator{\sdeg}{sdeg}

\DeclareMathOperator{\tor}{Tor}

\DeclareMathOperator{\inn}{in}

\newcommand{\ffi}{\varphi}

\newcommand{\fall}{\mbox{for all} ~}

\newcommand{\s}{\; | \;}

\def\pnt{{\raise0.5mm\hbox{\large\bf.}}}
\def\lpnt{{\hbox{\large\bf.}}}

\title[Algebraic shifting and exterior and symmetric algebra methods]
{Algebraic shifting and exterior and symmetric algebra methods}

\author{Uwe Nagel}
\address
{Department of Mathematics, University of Kentucky,
715 Patterson Office Tower, Lexington, KY 40506-0027, USA}
\email{uwenagel@ms.uky.edu}

\author{Tim R\"omer}
\address{FB Mathematik/Informatik, Universit\"at Osnabr\"uck, 49069
  Osnabr\"uck, Germany}
\email{troemer@mathematik.uni-osnabrueck.de}

\author{Natale Paolo Vinai$^*$}
\address{Department of Mathematics, University of Bologna, Piazza di
  Porta S. Donato, 5, 40126 Bologna,
Italy}
\email{vinai@dm.unibo.it}

\thanks{$^*$ The third author gratefully acknowledges partial support by
the
FIRB Research Project "Teoria dell'intersezione e applicazioni computazionali".
}

\begin{document}
\begin{abstract}
We define and study Cartan--Betti numbers of a graded ideal $J$ in the exterior algebra
over an infinite field which include the usual graded Betti numbers of $J$ as a special case.
Following ideas of Conca regarding Koszul--Betti numbers over the symmetric algebra,
we show that Cartan--Betti numbers increase by passing to the generic initial ideal and
the squarefree lexsegement ideal respectively.
Moreover, we characterize the cases where the inequalities become equalities.
As combinatorial applications
of the first part of this note and some further symmetric algebra methods
we establish results about algebraic shifting of simplicial complexes
and use them to compare different shifting operations.  In particular,
we show that each shifting operation does not decrease the number of
facets, and that the exterior shifting is the best among the exterior
shifting operations in the sense that it increases the number of
facets the least.
\end{abstract}

\maketitle
\section{Introduction}
Let $S=K[x_1,\dots,x_n]$ be a polynomial ring over a field $K$ of characteristic $0$,
let $I \subset S$ be a graded ideal and denote by
$\beta^S_{ij}(S/I)=\dim_K \Tor_i^S(S/I,K)_j$
the graded Betti numbers of $S/I$.
In the last decades
the graded Betti numbers $\beta^S_{ij}(S/I)$ were studied intensively.
To $I$ one associates several important
monomial ideals like the generic initial ideal $\gin(I)$ with respect to the
reverse lexicographic order, or the lexsegment ideal $\lex(I)$.
It is well-known by work of \cite{BI93}, \cite{HU93} and \cite{PA96}
that there are inequalities
$\beta^S_{ij}(S/I) \leq \beta^S_{ij}(S/\gin(I)) \leq \beta^S_{ij}(S/\lex(I))$ for all $i,j$.
The cases where these are equalities for all $i,j$ have been characterized by
Aramova--Herzog--Hibi \cite{ARHEHI00ideal} and \cite{HEHI99}.
The first inequality is an equality for all $i,j$ if and only if $J$ is a componentwise linear ideal,
and we have equalities everywhere for all $i,j$ if and only if $J$ is a Gotzmann ideal.
All these results were generalized to so--called Koszul--Betti numbers in \cite{CO04}.

Let $E=K\langle e_1,\dots,e_n \rangle$ be the
exterior algebra over an infinite field $K$ and $J \subset E$ a graded ideal.
Aramova--Herzog--Hibi (\cite{ARHEHI97}) showed that the constructions and results mentioned above
hold similarly for $J$. More precisely, there exists a generic initial ideal $\gin(J)\subset E$
with respect to the reverse lexicographic order,
and the unique squarefree lexsegment ideal $\lex(J)\subset E$ with the same Hilbert function as $J$.
Let $\beta^E_{ij}(E/J)=\dim_K \Tor_i^E(E/J,K)_j$ be the graded Betti numbers of $E/J$.
Then we also have the inequalities
$\beta^E_{ij}(E/J) \leq \beta^E_{ij}(E/\gin(J)) \leq \beta^E_{ij}(E/\lex(J))$
for all integers $i,j$. The first inequality is again an equality for all $i,j$
if and only if $J$ is a componentwise linear ideal
as was observed in \cite{ARHEHI00ideal}.

The first part of this paper
is devoted to extend the latter results to statements about Cartan--Betti numbers,
similarly as Conca did in the symmetric algebra case for Koszul--Betti numbers.
Over the exterior algebra there exists
the construction of the Cartan complex and Cartan homology (see Section \ref{cartan} for details),
which behave in many ways like the Koszul complex and Koszul homology.
Following Conca's ideas we define the Cartan--Betti numbers as
$
\beta^E_{ijp}(E/J)= \dim_K H_i(f_1,\dots,f_p;E/J)_j
$
where $H_i(f_1,\dots,f_p;E/J)$ is the $i$-th Cartan homology
of $E/J$ with respect to a generic sequence of linear forms
$f_1,\dots,f_p$ for $1\leq p \leq n$. Observe that this definition does not depend on the chosen
generic sequence and that these modules are naturally graded.
We set $H_i(f_1,\dots,f_p;E/J)=0$ for $i>p$. Note that
$\beta^E_{ijn}(E/J)= \beta^E_{ij}(E/J)$ are the graded Betti numbers of $E/J$.
Observe that a generic initial ideal $\gin_\tau(J)$ can be defined with respect to any term
order
$\tau$ on $E$.
Our first main result, Theorem \ref{mainleq}, shows that,  for all integers $i,j,p$:
$$
\beta^E_{ijp}(E/J) \leq \beta^E_{ijp}(E/\gin_\tau(J)) \leq
\beta^E_{ijp}(E/\lex(J)) .
$$
Analogously to result in \cite{CO04} we can characterize precisely
the cases where the inequalities become equalities.
At first we have:
\begin{thm}
Let $J \subset E$ be a graded ideal.
The following conditions are equivalent:
\begin{enumerate}
\item
$\beta_{ijp}^E(E/J)= \beta_{ijp}^E(E/\gin(J))$ for all $i,j,p$;
\item
$\beta_{1jn}^E(E/J)= \beta_{1jn}^E(E/\gin(J))$ for all $j$;
\item
$J$ is a componentwise linear ideal;
\item
A generic sequence of linear forms $f_1,\dots,f_n$ is a proper sequence of $E/J$.
\end{enumerate}
\end{thm}
We refer to Section \ref{cartan} for the definition of a proper sequence.
Moreover, we show:
\begin{thm}
For each graded ideal $J \subset E$, the following statements are
equivalent:
\begin{enumerate}
\item
$\beta^E_{ijp}(E/J)= \beta^E_{ijp}(E/\lex(J))$
for all $i,j,p$;
\item
$\beta^E_{1jn}(E/J)= \beta^E_{1jn}(E/\lex(J))$
for all $j$;
\item
$J$ is a Gotzmann ideal in the exterior algebra;
\item
$\beta^E_{0jp}(E/J) = \beta^E_{0jp}(E/\lex(J))$
for all $j,p$
and
$J$ is componentwise linear.
\end{enumerate}
\end{thm}
Section \ref{cartan} ends with the discussion which generic initial
ideal increases the Cartan--Betti numbers the least.
In Theorem \ref{mainginleq} we show that for all $i,j,p$ and for any term order $\tau$:
$$
\beta_{ijp}^E(E/J)\leq  \beta_{ijp}^E(E/\gin(J)) \leq
\beta_{ijp}^E(E/\gin_\tau(J)) .
$$
The second part of this note presents combinatorial applications
of the results mentioned above and some symmetric algebra methods.
More precisely,
we study properties of simplicial complexes under
shifting operations. Let $\mathcal{C}([n])$ be
the set of simplicial complexes on $n$ vertices.   Following
\cite{KA84}, \cite{KA01}, a shifting operation is a map
$\shift \colon \mathcal{C}([n]) \to \mathcal{C}([n])$ satisfying
certain conditions.  We refer to Section \ref{sec-def} for precise
definitions. Intuitively, shifting replaces each complex $\Gamma$ with
a combinatorially simpler complex $\Delta (\Gamma)$ that still
captures some properties of $\Gamma$. Shifting has become an important
technique that has been successfully applied in various contexts and
deserves to be investigated in its own right (cf., for example,
\cite{ARHEHI00}, \cite{ARHE}, \cite{BNT}, \cite{BNT2},
\cite{MH1}, \cite{MH2}, \cite{Nevo}).

Several shifting operations can
be interpreted algebraically. Denote by $I_{\Gamma}$ the
Stanley--Reisner ideal of $\Gamma$ in the polynomial ring $S =
K[x_1,\ldots,x_n]$. Then
symmetric shifting $\Delta^s$ can be realized by passing to a kind of
polarization of the generic initial ideal of $I_{\Gamma}$ with respect
to the reverse lexicographic order (cf.\ Example
\ref{exsymmetric}). Similarly, denote by $J_{\Gamma}$ the exterior
Stanley--Reisner ideal in the exterior algebra $E = K\langle
e_1,\dots,e_n\rangle$. The passage from $J_{\Gamma}$ to its generic
initial ideal with respect to any term order $\tau$ on $E$ leads to
the exterior algebraic shifting operation $\Delta^{\tau}$ (Example
\ref{exexterior}). If $\tau$ is the reverse-lexicographic order, we
call $\Delta^e := \Delta^{\tau}$ the exterior algebraic shifting.

While each shifting operation preserves the $f$-vector, the number of
facets may change. However, we show in Theorem \ref{thm-adeg-incr}
that it cannot decrease.
After recalling basic definitions in Section \ref{sec-def},
we first study symmetric algebraic shifting. To this end we use degree
functions. In Section \ref{sec-degree} we show that the smallest
extended degree is preserved under symmetric algebraic shifting whereas the
arithmetic degree $\adeg (\Gamma)$ may increase because it equals the
number of facets. However, the smallest extended degree and the
arithmetic degree of each shifted complex agree.
Furthermore we show that exterior algebraic shifting shift is the best exterior shifting operation in the sense
that it increases the number of facets the least, i.\ e.\ for any term order $\tau$
$$
\adeg \Gamma \leq \adeg \Delta^e(\Gamma) \leq \adeg \Delta^\tau(\Gamma).
$$
Moreover, $\Delta^e$ preserves the arithmetic degree of $\Gamma$ if
and only $\Gamma$ is sequentially Cohen--Macaulay. These results rely
on the full strength of the results in Section \ref{cartan}. The proof
also uses Alexander duality and a reinterpretation of the arithmetic
degree over the exterior algebra.

In Section \ref{ringprop} we use degree functions to give a
short new proof of the fact that $\Delta^s (\Gamma)$ is a pure complex
if and only if $\Gamma$ is Cohen--Macaulay. We then derive a
combinatorial interpretation of the smallest extended degree using
iterated Betti numbers and show that these numbers agree with the
$h$-triangle if and only if $\Gamma$ is sequentially Cohen--Macaulay.


For notions and results related to commutative algebra we refer to \cite{EI}, \cite{BRHE98} and \cite{VA}.
For details on
combinatorics we refer to the book \cite{ST96}
and \cite{KA01}.

\section{Cartan--Betti numbers}
\label{cartan}
In this section we introduce and study Cartan--Betti numbers of
graded ideals in the exterior algebra $E=K\langle e_1,\dots,e_n \rangle$ over
an infinite field $K$. Our results rely on techniques from Gr\"obner basis
theory in the exterior algebra. Its basics are treated in \cite{ARHEHI97}. We begin with
establishing some extensions of the theory that are analogous to results over the
polynomial ring.

Recall that a monomial of degree $k$ in $E$ is an element
$e_F=e_{a_1} \wedge \ldots \wedge e_{a_k}$
where $F=\{a_1,\ldots,a_k\}$ is a subset of $[n]$ with $a_1< \ldots < a_k$.
To simplify notation, we write sometimes $f g = f \wedge g$ for any two elements $f, g \in E$.
We will use only term orders $\tau$ on $E$ that satisfy
$
e_1 >_{\tau} e_2 >_{\tau} \ldots >_{\tau} e_n.$
Given a term order $\tau$ on $E$
and a graded ideal $J\subset E$
we denote by $\ini_\tau(J)$ and $\gin_\tau(J)$ respectively
the {\em initial ideal} of $J$ and the {\em generic initial ideal} of $J$
with respect to $\tau$.
We also write
$\ini(J)$ and $\gin(J)$ for
the initial ideal of $J$ and the generic initial ideal of $J$
with respect to the reverse lexicographic order on $E$ induced by
$e_1>\cdots>e_n$.

Recall that a monomial ideal
$J \subset E$ is called a {\em squarefree strongly stable ideal} with respect to $e_1>\dots>e_n$
if for all $F\subseteq [n]$ with $e_F \in J$
and all $i \in F$, $j<i$, $j \not\in F$
we have $e_j\wedge e_{F \setminus\{i\}} \in J$.
Given an ideal $J\subset E$ and a term order $\tau$ on $E$, note that
$\gin_\tau(J)$ is
always squarefree strongly stable and that $\gin_\tau(J)=J$ iff $J$ is
squarefree
strongly stable which can be seen analogously to the case of a
polynomial ring. As
expected, the passage to initial ideals increases the graded Betti
numbers. The next
result generalizes \cite[Proposition 1.8]{ARHEHI97}. It is the
exterior analogue of a well-known fact  for the polynomial ring (cf.\
\cite[Lemma 2.1]{CO04}).

\begin{prop}
\label{helpergb}
Let $J, J' \subset E$ be graded ideals and let $\tau$ be any term
order on $E$. Then, for all $i, j \in \Z$:
$$
\dim_K [\tor^E_i (E/J, E/J')]_j \leq \dim_K [\tor^E_i (E/\inn_{\tau}
  (J), E/\inn_{\tau} (J'))]_j.
$$
\end{prop}

\begin{proof}
This follows from standard deformation arguments. The proof is
verbatim the same as the
one of the symmetric version in \cite[Lemma 2.1]{CO04} after replacing
the polynomial
ring $S$ with the exterior algebra $E$.
\end{proof}

In the following we present exterior versions of results of Conca \cite{CO04}.
We follow the strategy of his proofs and replace symmetric algebra methods and the Koszul
complex by exterior algebra methods and the Cartan complex.  We first recall the
construction of the Cartan complex (see, e.g., \cite{ARHEHI97} or \cite{HE01}). For
important results on this complex, we refer to \cite{ARAVHE00}, \cite{ARHEHI97}, \cite{ARHEHI98},
\cite{ARHEHI00ideal}.

For a sequence $\vb=v_1,\ldots, v_m \subseteq E_1$,  the Cartan complex
$C_{\lpnt}(\vb;E)$ is defined to be the free divided power algebra $E\langle
x_1,\ldots,x_m \rangle$ together with a differential $\delta$. The free divided power
algebra $E\langle x_1,\ldots,x_m \rangle$ is generated over $E$ by the divided powers
$x_{i}^{(j)}$, \ $i=1,\ldots,m$ and $j\geq 0$, that satisfy the relations
$x_{i}^{(j)}x_{i}^{(k)}=\frac{(j+k)!}{j!k!}x_{i}^{(j+k)}$. We set $x_{i}^{(0)}=1$ and
$x_{i}^{(1)}=x_{i}$ for $i=1,\ldots, m$. Hence $C_{\lpnt}(\vb;E)$ is a free $E$-module
with basis $x^{(a)}:=x_{1}^{(a_1)}\ldots x_{m}^{(a_m)}$, $a = (a_1,\ldots,a_m) \in
\N^{m}$. We set $\deg x^{(a)}=i$ if $|a|:=a_1+\ldots+a_m=i$ and
$C_{i}(\vb;E)=\bigoplus_{|a|=i} Ex^{(a)}$. The $E$-linear differential $\partial$ on
$C_{\lpnt}(\vb;E)$ is defined as follows. For $x^{(a)}=x_{1}^{(a_1)}\ldots x_{m}^{(a_m)}$
we set $\partial(x^{(a)})=\sum_{a_i>0} v_{i} \cdot x_{1}^{(a_1)}\ldots x_{i}^{(a_{i}-
1)}\ldots x_{m}^{(a_m)}.$ One easily checks  that $\partial\circ\partial=0$, thus
$C_{\lpnt}(\vb;E)$ is indeed a complex.

We denote by $\mathcal{M}$ the category of finitely generated graded left and right
$E$-modules $M$, satisfying $ax=(-1)^{|\deg(a)||\deg(x)|}xa$ for all homogeneous $a\in E$
and $x\in M$. For example, every graded ideal $J\subseteq E$ belongs to $\mathcal{M}$.
For $M \in \mathcal{M}$ and $\vb=v_1,\ldots, v_m \subseteq E_1$, the complex
$C_{\lpnt}(\vb;M)=C_{\lpnt}(\vb;E)\otimes_E M$ is called the {\em Cartan complex} of $M$
with respect to $\vb$. Its homology is denoted by $H_{\lpnt}(\vb;M)$; it is called {\em
Cartan homology}.

There is a natural grading of this complex and its homology. We set $\deg x_i =1$ and
$C_{j}(\vb;M)_{i}:=\text{span}_{K}\{m_{a}x^{(b)} \s m_a \in M_a, a+|b|=i,|b|=j \}.$
Cartan homology can be computed recursively. For $j=1,\ldots,m-1$, the following sequence
is exact:
\begin{eqnarray} \label{eq-cartan}
\hspace*{.5cm} 0 \to C_{\lpnt}(v_1,\ldots,v_j;M) \overset{\iota}{\to}
C_{\lpnt}(v_1,\ldots,v_{j+1};M)\overset{\ffi}{\to}C_{\lpnt-
1}(v_1,\ldots,v_{j+1};M)(-1)\to 0.
\end{eqnarray}
Here $\iota$ is a natural inclusion map and $\ffi$ is given by
$$
\ffi(g_0 + g_1x_{j+1}+\ldots+g_{k}x_{j+1}^{(k)})=g_1 + g_2x_{j+1}+\ldots+g_{k}x_{j+1}^{(k-1)}
\text{ for }
g_i \in C_{i-k}(v_1,\ldots ,v_{j};M).
$$
The associated long exact homology sequence is
$$
\ldots \to H_{i}(v_1,\ldots,v_j;M) \overset{\alpha_i}{\to}
H_{i}(v_1,\ldots,v_{j+1};M)\overset{\beta_i}{\to} H_{i-
1}(v_1,\ldots,v_{j+1};M)(-1)
$$
$$
\overset{\delta_{i-1}}{\to} H_{i-1}(v_1,\ldots,v_j;M) \overset{\alpha_{i-
1}}{\to} H_{i-1}(v_1,\ldots,v_{j+1};M)\overset{\beta_{i-1}}{\to}\ldots,
$$
where $\alpha_{i}$ and $\beta_{i}$ are induced by $\iota$ and $\ffi$, respectively. For a
cycle $z=g_0 + g_1x_{j+1}+\ldots+g_{i-1}x_{j+1}^{(i-1)}$ in $C_{i-
1}(v_1,\ldots,v_{j+1};M)$ one has $\delta_{i-1}([z])=[g_0v_{j+1}]$.
It is now easy to
see that e.~ g.\ for
$e_t,\dots,e_n$, the Cartan complex $C_{\lpnt}(e_t,\dots,e_n;E)$ is a free
resolution of $E/(e_t,\dots,e_n)$. Hence, for each module $M \in
\mathcal{M}$, there are isomorphisms
\begin{equation}
\label{tor}
\Tor_i^{E}(E/(e_t,\dots,e_n),M)\cong H_{i}(e_t,\dots,e_n;M).
\end{equation}
In particular,
for $\eb=e_1,\dots,e_n$ we have $K\cong E/(\eb)$ and there are isomorphisms of graded $K$-vector spaces $ \Tor_i^{E}(K,M)\cong
H_{i}(\eb;M)$.

Cartan homology is useful to study resolutions of graded ideals in the exterior algebra.
For example, it has been used to show that the Castelnuovo--Mumford regularity $\reg(E/J)=\max\{j \mid \exists $ $i \text{ such that } \Tor_i(K,E/J)_{i+j}\neq 0 \}$
does not change by passing to the generic initial ideal with respect to the reverse
lexicographic order, i.~e.\ $ \reg(E/J)= \reg(E/\gin(J))$ (cf.\ \cite{ARHE}, Theorem 5.3).

There are several other algebraic invariants which behave similarly. See for example
\cite{HETE} for results in this direction. Here, we are interested in the following
numbers which, for $p = n$, include the graded Betti numbers of $E/J$:

\begin{defi} \label{def-CB}
Let $J \subset E$ be a graded ideal and
$f_1,\dots,f_p$ be a generic sequence of linear forms.
Let $C(f_1,\dots,f_p;E/J)$ be the Cartan complex
of $E/J$ with respect to that sequence
and $H(f_1,\dots,f_p;E/J)$ the corresponding Cartan homology.
We denote by
$$
\beta^E_{ijp}(E/J)= \dim_K H_i(f_1,\dots,f_p;E/J)_j
$$
the {\em Cartan--Betti numbers} of $E/J$ for $i=0,\dots,p$.
We set $H_i(f_1,\dots,f_p;E/J)=0$ for $i>p$.
\end{defi}

A generic sequence means here that
there exists a non-empty Zariski open set $U$ in $K^{p\times n}$ such that
if one chooses the $p\times n$ coefficients
that express  $f_1,\dots,f_p$ as linear combinations of  $e_1,\dots,e_n$ inside $U$,
then the ranks of the maps in each degree in  the Cartan complex are as big as possible.
In particular, the $K$-vector space dimension of the  homology
is the same for each generic sequences and thus the definition of
$\beta^E_{ijp}(E/J)$ does not depend on the chosen generic sequence.

We are going to explain that it is easy to compute the Cartan--Betti numbers for
squarefree strongly stable ideals in the exterior algebra. In fact, given a generic
sequence of linear forms $f_1,\dots,f_p$ there exists an invertible upper triangular
matrix $g$ such that the linear space spanned by $f_1,\dots,f_p$ is mapped by the induced
isomorphism $g \colon E \to E $ to the one generated by
$e_{n-p+1},\dots,e_n$. It follows
that, as graded
$K$-vector spaces, $ H_i(f_{1},\dots,f_p;E/J)\cong
H_i(e_{n-p+1},\dots,e_n;E/g(J))$. Since $J$ is squarefree strongly
stable we have that $g(J)=J$ and thus we get:

\begin{lem}
\label{cartanhelper} Let $J \subset E$ be a squarefree strongly stable ideal. Then we
have for all  $i,j,p$:
$$
\beta_{ijp}^E(E/J)=\dim_K H_i(e_{n-p+1},\dots,e_n;E/J)_j.
$$
\end{lem}

Aramova, Herzog and Hibi (\cite{ARHEHI97}) computed the Cartan--Betti numbers for squarefree strongly
stable ideals which gives rise to Eliahou--Kervaire type resolutions in the exterior
algebra. To present this result we need some more notation. Recall that for a monomial
ideal $J\subset E$ we denote by $G(J)$ the unique minimal system of monomial generators
of $J$. For a monomial $e_F \in E$ where $F \subseteq [n]$ we set
$\max(F)=\max\{i \in F\}.$ Given a set $G$ of monomials we define:
\begin{eqnarray*}
m_i(G) & = & |\{e_F \in G : \max(e_F)=i  \}|, \\
m_{\leq i}(G) &  = & |\{e_F \in G : \max(e_F) \leq i \}|, \\
m_{ij}(G) & = & |\{e_F \in G : \max(e_F)=i,\ |F|=j  \}|.
\end{eqnarray*}
For a monomial ideal or a vector space generated by monomials $J$, we denote by $m_i(J)$,
$m_{\leq i}(J)$ and $m_{ij}(J)$ the numbers $m_i(G)$, $m_{\leq i}(G)$ and $m_{ij}(G)$
where $G$ is the set of minimal monomial generators of $J$. Using Lemma
\ref{cartanhelper} we restate the result of Aramova, Herzog and Hibi as follows:

\begin{thm}
\label{cartanbetti}
Let $J \subset E$ be a squarefree strongly stable ideal.
Then the Cartan--Betti numbers are given by
the following formulas:
\begin{enumerate}
\item{(Aramova--Herzog--Hibi)}
For all $i>0$, $p \in [n]$ and every $j \in \Z$
we have
$$
\beta_{ijp}^E(E/J)
=
\sum_{k = n-p+1}^n m_{k,j-i+1}(J) \binom{k+p-n+i-2}{i-1}.
$$
\item
For all $p \in [n]$ and every $j \in \Z$
we have
$$
\beta_{0jp}^E(E/J)
=
\binom{n-p}{j}-m_{\leq n-p}(J_j)
$$
where $J_j$ is the degree $j$ component of the ideal $J$.
\end{enumerate}
\end{thm}
\begin{proof}
Only (ii) needs a proof since (i) immediately follows from the proof of Proposition 3.1 of \cite{ARHEHI97}. But since $e_{n-p+1},\dots,e_n$ is generic for $E/J$
we can compute  $\beta_{0jp}^E(E/J)$
from the Cartan homology with respect to this sequence. Then the formula
follows from a direct computation.
\end{proof}

Note that for $p=n$ we get the graded Betti numbers of $E/J$. 
Using the numbers $m_{\leq i}(J_j)$, it is possible to compare the Cartan--Betti numbers
of squarefree strongly stable ideals. In fact we have the following:

\begin{prop}
\label{maintechnical}
Let $J, J' \subset E$
be squarefree strongly stable ideals with the same Hilbert function
such that
$m_{\leq i}(J_j)\leq m_{\leq i}(J'_j)$ for all $i,j$.
Then
$
\beta_{ijp}^E(E/J')
\leq
\beta_{ijp}^E(E/J)
$
for all
$i,j,p$
with equalities everywhere
if and only if
$m_{\leq i}(J_j)= m_{\leq i}(J'_j)$ for all $i,j$.
\end{prop}

\begin{proof}
It follows from
Theorem \ref{cartanbetti} (ii) and the assumption that
\begin{eqnarray*}
\beta_{0jp}^E(E/J) - \beta_{0jp}^E(E/J') &=& \binom{n-p}{j} - m_{\leq n-p}(J_j) -
\binom{n-p}{j} + m_{\leq n-p}(J'_j)\\
&=&
m_{\leq n-p}(J'_j) - m_{\leq n-p}(J_j)\\
&\geq&
0.
\end{eqnarray*}
This shows the desired inequalities for $i=0$.
Next assume that $i>0$ and observe 
\begin{eqnarray*}
\lefteqn{ m_{k}(J_{j-i+1})=|\{e_F \in J : \max (e_F) = k, \  |F| = j-i+1\}| }\\
&=&
|\{e_F \in G(J) : \max (e_F) = k,  |F| = j-i+1\}| \\
&+&
| \{e_F \in J : \exists e_H \in J_{j-i},  \max (e_H)\leq k-1, e_F=e_H \wedge e_k\} |\\
&=&
m_{k,j-i+1}(J) + m_{\leq k-1}(J_{j-i}).\\
\end{eqnarray*}

We compute

\begin{eqnarray*}
\lefteqn{ \beta_{ijp}^E(E/J') } \\
&=&
\sum_{k = n-p+1}^n m_{k,j-i+1}(J') \binom{k+p-n+i-2}{i-1}\\
&=&
\sum_{k = n-p+1}^n \bigl( m_{k}(J'_{j-i+1}) - m_{\leq k-1}(J'_{j-i}) \bigr) \binom{k+p-n+i-2}{i-1}\\
&=&
\sum_{k = n-p+1}^n \bigl(
m_{\leq k}(J'_{j-i+1})- m_{\leq k-1}(J'_{j-i+1}) - m_{\leq k-1}(J'_{j-i})
\bigr) \binom{k+p-n+i-2}{i-1}\\
&=&
m_{\leq n}(J'_{j-i+1}) \binom{p+i-2}{i-1}
- m_{\leq n-p}(J'_{j-i+1})\\
&+&
\sum_{k = n-p+1}^{n-1}
m_{\leq k}(J'_{j-i+1})
\bigl(
\binom{k+p-n+i-2}{i-1}
-
\binom{k+p-n+i-1}{i-1}
\bigr)\\
&-&
\sum_{k = n-p+1}^n
m_{\leq k-1}(J'_{j-i})
\binom{k+p-n+i-2}{i-1}\\
&=&
\dim_K (J'_{j-i+1}) \binom{p+i-2}{i-1}
- m_{\leq n-p}(J'_{j-i+1})\\
&-&
\sum_{k = n-p+1}^{n-1}
m_{\leq k}(J'_{j-i+1})
\binom{k+p-n+i-2}{i-2}
\\
&-&
\sum_{k = n-p+1}^n
m_{\leq k-1}(J'_{j-i})
\binom{k+p-n+i-2}{i-1}\\
\end{eqnarray*}
and similarly $\beta_{ijp}^E(E/J)$. We know by assumption that $\dim_K
(J'_{j-i+1})=\dim_K (J_{j-i+1})$ and $m_{\leq i}(J_j)\leq m_{\leq i}(J'_j)$ for all
$i,j$. Hence  the last equation of the computation implies  $\beta_{ijp}^E(E/J)-
\beta_{ijp}^E(E/J') \geq 0$ because the left-hand side is a sum of various $m_{\leq
i}(J'_j)-m_{\leq i}(J_j)$ times some binomial coefficients. Moreover, we have equalities
everywhere if and only if $m_{\leq i}(J_j)= m_{\leq i}(J'_j)$ for all $i,j$.
\end{proof}

Using the above results as well as  \cite[Lemma 3.7]{ARHEHI98} instead of Proposition 3.5
in \cite{CO04}, we get the following rigidity statement that is proved analogously to
Proposition 3.7 in \cite{CO04}:

\begin{cor}
\label{rigid}
Let $J, J' \subset E$
be a squarefree strongly stable ideals with the same Hilbert function
and $m_{\leq i}(J_j)\leq m_{\leq i}(J'_j)$ for all $i,j$.
Then the following statements are equivalent:
\begin{enumerate}
\item
$\beta_{ijp}^E(E/J')=\beta_{ijp}^E(E/J)$
for all $i,j$ and $p$;
\item
$\beta_{ij}^E(E/J')=\beta_{ij}^E(E/J)$
for all $i$ and $j$;
\item
$\beta_{1j}^E(E/J')=\beta_{1j}^E(E/J)$
for all $j$;
\item
$\beta_{1}^E(E/J')=\beta_{1}^E(E/J)$;
\item
$m_{i}(J'_j)=m_{i}(J_j)$
for all $i,j$;
\item
$m_{\leq i}(J'_j)=m_{\leq i}(J_j)$
for all $i,j$.
\end{enumerate}
\end{cor}

Given a graded ideal $J \subset E$ there exists the unique
{\em squarefree lexsegment ideal} $\lex(J) \subset E$ with the same Hilbert function as $J$
(see \cite{ARHEHI97} for details).
Next we present a mild variation of a crucial result
in \cite{ARHEHI98} which was an important step to compare Betti numbers
of squarefree strongly stable ideals
and their corresponding squarefree lexsegment ideals.

\begin{thm}{(Aramova--Herzog--Hibi)}
\label{lexext}
Let $J \subset E$ be a squarefree strongly stable ideal and let $L \subset
E$ be its squarefree lexsegment ideal. Then
$$m_{\leq i}(L_j) \leq m_{\leq i}(J_j) \text{ for all } i,j.$$

\end{thm}
\begin{proof}
It is a consequence of \cite[Theorem 3.9]{ARHEHI98} by observing the following facts.
For $j\in \N$ consider the ideals $J'=(J_j)$ and $L'=(L_j)$ of $E$.
Then $J'$ is  a squarefree strongly stable ideal,
$L'$ is a squarefree lexsegment ideal,
and
$\dim_K L'_t \leq \dim_K J'_t$ for all $t$ by \cite[Theorem 4.2]{ARHEHI97}.
\end{proof}

As \cite[Lemma 2]{CO03} one shows for Cartan--Betti numbers:

\begin{lem}
\label{someeq}
Let $J \subset E$ be a graded ideal.
Then
$$
\beta_{0jp}^E(E/J) = \beta_{0jp}^E(E/\gin(J)) \text{ for all } j,p.
$$
\end{lem}

\begin{proof}
Consider the graded  $K$-algebra homomorphism $g\colon E\to E$ that is induced by a
generic matrix in $\GL_n(K)$. Then the linear forms
$f_i := g^{-1}(e_{n-p+i})$, $i=1,\ldots,p$, are generic too and the Hilbert functions of
$E/(J+L)$ and $E/(g(J)+L')$
coincide, where $L=(f_1,\dots,f_p)$ and $L' = (e_{n-p+1},\ldots,e_n)$. Note that passing
to the initial ideals does not change the Hilbert function and that
$\ini(g(J)+L')=\ini(g(J))+L'$
because the chosen term order is revlex. Furthermore, we
have that $\ini(g(J))=\gin(J)$ because $g$ is generic. Hence, the Hilbert function of
$E/(J+L)$
equals to that of
$E/(\gin(J)+L')$.
Since
$ \beta_{0jp}^E(E/J)= \dim_K (E/(J+L))_j $
and
$ \beta_{0jp}^E(E/\gin(J))= \dim_K (E/(\gin(J)+L'))_j $,
our claim follows.
\end{proof}

The following result generalizes \cite[Prop. 1.8, Theorem 4.4]{ARHEHI97}. The proof is
similar to \cite[Theorem 4.3]{CO04}.

\begin{thm}
\label{mainleq}
Let $J \subset E$ be a graded ideal and
$\tau$ an arbitrary term order on $E$.
Then
\begin{enumerate}
\item
$
\beta_{ijp}^E(E/J)
\leq
\beta_{ijp}^E(E/\gin_\tau(J))
$
for all $i,j,p$.
\item
$
\beta_{ijp}^E(E/J)
\leq
\beta_{ijp}^E(E/\lex(J))
$
for all $i,j,p$.
\end{enumerate}
\end{thm}
\begin{proof}
(i):  Let $g\in \GL_n(K)$ be a generic matrix and let $f_i$ be the preimage of
$e_{n-p+1}$ under the induced $K$-algebra isomorphism $g \colon E\to E$ where
$i=1,\dots,p$. Then
$$
H_i(f_1,\dots, f_p;E/J)
\cong
H_i(e_{n-p+1},\dots,e_{n};E/g(J))
\cong
\Tor^E_i(E/g(J), E/(e_{n-p+1},\dots,e_{n}))
$$
where the last isomorphism was noted above in (\ref{tor}).
Proposition \ref{helpergb} provides:
\begin{eqnarray*}
\beta^E_{ijp}(E/J)
&\leq&
\dim_K \Tor^E_i(E/\ini_\tau(g(J)), E/(e_{n-p+1},\dots,e_{n}))_j\\
&=&
\dim_K  H_i(e_{n-p+1},\dots,e_{n}; E/\ini_\tau(g(J)))_j.
\end{eqnarray*}
Since $g$ is a generic matrix, $\ini_\tau(g(J))=\gin_\tau(J)$  is
squarefree strongly
stable. Thus we can apply Lemma \ref{cartanhelper} and (i) follows.

(ii): By  (i) we may replace $J$ by $\gin(J)$, thus we may assume that
$J$ is squarefree
strongly stable. Now (ii) follows from  \ref{lexext} (ii) and \ref{maintechnical}.
\end{proof}

The next results answer the natural question, in which cases we have
equalities in Theorem \ref{mainleq}. Denote by $J_{\langle t\rangle}$ the ideal that is generated by all degree
$t$ elements of the graded ideal $J$. Then  $J \subset E$ is called 
{\em componentwise linear} if, for all $t\in \N$, the ideal  has an $t$-linear resolution, 
i.~e.\ $\Tor_i^E(J_{\langle t\rangle},K)_{i+j}=0$ for $j \neq t$.

Given a module $M \in \mathcal{M}$, we call, in analogy to the symmetric case, 
the sequence $v_1,\dots,v_m \in E$ a {\em proper $M$-sequence} if
for all $i\geq 1$ and all $1\leq j< m$ the maps
$$ 
\delta_i \colon H_i(v_1,\dots,v_{j+1};M)(-1) \to H_i(v_1,\dots,v_{j};M)
$$
are zero maps.
The following result generalizes \cite[Theorem 2.1]{ARHEHI00ideal}. 
We follow the idea of the proof of \cite[Theorem 4.5]{CO04}.

\begin{thm}
\label{maincl}
Let $J \subset E$ be a graded ideal.
The following conditions are equivalent:
\begin{enumerate}
\item
$\beta_{ijp}^E(E/J)= \beta_{ijp}^E(E/\gin(J))$ for all $i,j,p$;
\item
$\beta_{1jn}^E(E/J)= \beta_{1jn}^E(E/\gin(J))$ for all $j$;
\item
$J$ is a componentwise linear ideal;
\item
A generic sequence of linear forms $f_1,\dots,f_n$ is a proper sequence of $E/J$.
\end{enumerate}
\end{thm}
\begin{proof}
(i) $\Rightarrow$ (ii): This is trivial.

(ii) $\Rightarrow$ (iii): This follows from of \cite[Theorem 1.1 and Theorem
2.1]{ARHEHI00ideal}. In fact, the proof of Aramova, Herzog an Hibi shows  this stronger
result.

(iii) $\Rightarrow$ (iv):
Assume that $J$ is componentwise linear.
Let $f_1,\dots,f_n$ be a generic sequence of linear forms.
We have to show that 
for all $i\geq 1$ and all $1\leq j< n$
the homomorphisms
$ 
\delta_i \colon H_i(f_1,\dots,f_{j+1};E/J)(-1) \to H_i(f_1,\dots,f_{j};E/J)$
are zero maps.

Assume first that $J$ is generated in a single degree $d$. Since the regularity does not
change by passing to the generic initial ideal with respect to revlex, it follows that
also $\gin(J)$ is generated in degree $d$. By Theorem \ref{cartanbetti} (i) we get that
for $i\geq 1$ the homology module $H_i(f_1,\dots,f_j;E/\gin(J))$ is concentrated in the
single degree $i+d-1$. Theorem \ref{mainleq} (i) shows that the same is  true also for $H_i(f_1,\dots,f_j;E/J)$.
This implies that $H_i(f_1,\dots,f_{j+1};E/J)(-1)$ has non-trivial homogeneous elements only in degree $i+d$
and the module
$H_i(f_1,\dots,f_{j};E/J)$ has only homogeneous elements in degree $i+d-1$.
Since $\delta_i$ is a homogeneous homomorphism of degree zero, this implies that
simply by degree reasons $\delta_i$ is the zero map.


Now we consider the general case. Recall that
for a cycle $c$ in some complex we write $[c]$ for
the corresponding homology class.
For any element $e \in E$ we denote by $\bar e$ the corresponding
residue class in $E/J$ to distinguish notation.
Now we fix some $i$ and $j$. 
Let $a \in H_i(f_1,\dots,f_j;E/J)(-1)$ be a homogeneous element of degree $s$. 
We have to show that $\delta_i(a)=0$.
There exists an homogeneous element $z \in C_i(f_1,\dots,f_j;E)$ of degree $s-1$ 
such that $\bar z \in C_i(f_1,\dots,f_j;E/J)$ is a cycle representing $a$,
i.~e. $\partial_i(\bar z)=0$ and $a=[\bar z]$.
Since $\partial_i(\bar z)=0$ we have that
$\partial_i(z) \in JC_{i-1}(f_1,\dots,f_j;E)$ because 
$JC_{i-1}(f_1,\dots,f_j;E)=C_{i-1}(f_1,\dots,f_j;J)$
and we have the commutative diagram

$$
\begin{array}{ccccccc}
0  \ \rightarrow & C_i(f_1,\dots,f_j;J) & \longrightarrow & C_i(f_1,\dots,f_j;E) & \longrightarrow & C_i(f_1,\dots,f_j;E/J) & \rightarrow \ 0\\
                  & \downarrow \partial_i     &            &    \downarrow \partial_i &           & 
\downarrow \partial_i  &                   \\
0 \ \rightarrow & C_{i-1}(f_1,\dots,f_j;J) & \longrightarrow & C_{i-1}(f_1,\dots,f_j;E) & \longrightarrow & C_{i-1}(f_1,\dots,f_j;E/J) & \rightarrow \ 0\\
\end{array} 
$$
\noindent Since $\partial_i (z) $ is homogeneous of degree $s-1$, we obtain that $\partial_i(z)$ 
is in $J_k C_{i-1}(f_1,\dots,f_j;E)$ for $k=(s-1)-(i-1)=s-i$. 
Set $J' := J_{<k>}$. 
Then $[\bar z]$ can also be considered as the homology class  of $z$ in $H_i(f_1,\dots,f_j;E/J')$. 
By construction, $J'$ is generated in a single
degree and, by assumption, it has a linear resolution. Hence we know already that 
$\delta_i([\bar z])=0$ as an element of $H_i(f_1,\dots,f_{j-1};E/J')$ (in degree $s$). 
The natural homomorphism 
$H_i(f_1,\dots,f_{j-1};E/J') \to H_i(f_1,\dots,f_{j-1};E/J)$
induced by 
the short exact sequence $0\to J/J' \to E/J' \to E/J \to 0$
implies that
$\delta_i(a)=\delta_i([\bar z])=0$ as an element of $H_i(f_1,\dots,f_j;E/J)$ (in degree $s$) 
which we wanted to show.

(iv) $\Rightarrow$ (i):
By Lemma \ref{someeq} we know that
$\beta^E_{0jp}(E/J)=\beta^E_{0jp}(E/\gin(J))$ for all $j,p$.
Observe that
$\beta^E_{ij1}(E/J)=\beta^E_{ij1}(E/\gin(J))$ for $i>0$ and $j\in \Z$
as one easily checks using properties of $\gin$.
We prove (i) by showing
that the numbers $\beta_{ijp}^E(E/J)$ only depend on the numbers
$\beta_{0jp}^E(E/J)$ and $\beta_{ij1}^E(E/J)$.

Since $f_1,\dots,f_n$ is a proper sequence of $E/J$, the long exact Cartan homology
sequence splits into short exact sequences
$$
0 \to H_1(f_1,\dots,f_{p-1};E/J) \to H_1(f_1,\dots,f_{p};E/J) \to
H_0(f_1,\dots,f_{p};E/J)(-1)
$$
$$
\to
H_0(f_1,\dots,f_{p-1};E/J)
\to
H_0(f_1,\dots,f_{p};E/J)
\to
0
$$
and, for $i>1$,
$$
0 \to H_i(f_1,\dots,f_{p-1};E/J) \to H_i(f_1,\dots,f_{p};E/J) \to
H_{i-1}(f_1,\dots,f_{p};E/J)(-1) \to 0.
$$
Thus
$$
\beta^E_{1jp}(E/J) = \beta^E_{1jp-1}(E/J) + \beta^E_{0 j-1 p}(E/J) - \beta^E_{0jp-1}(E/J)
+ \beta^E_{0jp}(E/J)
$$
and, for $i > 1$,
$$
\beta^E_{ijp}(E/J) = \beta^E_{ijp-1}(E/J) + \beta^E_{i-1 j-1 p}(E/J).
$$
Since $\beta^E_{ijp}(E/J)=0$ for all $i>p$, it is easy to see that these equalities imply
(i).
\end{proof}

Recall that a graded ideal $J \subset E$ is called a {\em Gotzmann ideal} if the growth
of its Hilbert function is the least possible, i.\ e.\
$\dim_K E_1 \cdot J_j = \dim_K E_1\cdot  \lex(J)_j$
for all $j$. The next result is an exterior version of the corresponding result
\cite[Theorem 4.6]{CO04} in the polynomial ring.

\begin{thm}
\label{maingotz} For each graded ideal $J \subset E$, the following statements are
equivalent:
\begin{enumerate}
\item
$\beta^E_{ijp}(E/J)= \beta^E_{ijp}(E/\lex(J))$
for all $i,j,p$;
\item
$\beta^E_{1jn}(E/J)= \beta^E_{1jn}(E/\lex(J))$
for all $j$;
\item
$J$ is a Gotzmann ideal in the exterior algebra;
\item
$\beta^E_{0jp}(E/J) = \beta^E_{0jp}(E/\lex(J))$
for all $j,p$
and
$J$ is componentwise linear.
\end{enumerate}
\end{thm}

\begin{proof}

(i)  $\Rightarrow$  (ii): This is trivial.

(ii) $\Leftrightarrow$  (iii):
This follows immediately from the definition of Gotzmann ideals.

(ii) $\Rightarrow$  (iv): By Theorems \ref{mainleq} and \ref{maincl},  a Gotzmann ideal  is componentwise linear and
$\beta^E_{0jp}(E/J) = \beta^E_{0jp}(E/\lex(J))$ for all $j,p$.

(iv) $\Rightarrow$  (i):
From Theorem \ref{cartanbetti} (ii) we know that $m_{\leq n-p}(\gin(J)_j) = m_{\leq n-p}(\lex(J)_j)$ for all $j, p$. Then it follows from Corollary \ref{rigid} that
$\beta^E_{ijp}(E/\gin(J)) =  \beta^E_{ijp}(E/\lex(J))$  for all $i,j,p$. Since $J$ is componentwise linear, the assertion follows from Theorem \ref{maincl}.

\end{proof}

We conclude this section by comparing the Cartan--Betti numbers of
generic initial  ideals. The goal is to show that in the family
$$
\gins(J) :=\{\gin_\tau(J) : \tau \text{ a term order of } E\}
$$
of all generic initial ideals of $J$, the revlex-gin has the smallest
Cartan--Betti numbers.

We need a lemma and some more notation. Let $V \subset E_i$ be a
$d$-dimensional
subspace. Then $\bigwedge^d V$ is a $1$-dimensional subspace of
$\bigwedge^d E_i$. We
identify it with any of its nonzero elements. An exterior monomial in $\bigwedge^d E_i$
is by definition an element of the form $m_1 \wedge \ldots \wedge m_d$ where
$m_1,\ldots,m_d$ are distinct monomials in $E_i$. It is called a {\em $\tau$-standard}
exterior monomial if $m_1
>_{\tau} \ldots >_{\tau} m_d$. The vector space $\bigwedge^d E_i$ has a $K$-basis of
$\tau$-standard exterior monomials that we order lexicographically by
$
m_1 \wedge \ldots \wedge m_d >_{\tau} n_1 \wedge \ldots \wedge n_d
$
if $m_i >_{\tau} n_i$ for the smallest index $i$ such that $m_i \neq n_i$.
Using this order one defines the initial (exterior) monomial
$\inn_{\tau} (f)$ of any $f \in \bigwedge^d E_i$
as the maximal $\tau$-standard exterior monomial with respect to $>_{\tau}$
which
appears in the unique representation of $f$ as a sum of
$\tau$-standard exterior monomials.
Similarly, the initial subspace $\inn_{\tau}(V)$ of any subspace $V$ of $\bigwedge^d E_i$ is defined as the subspace
generated by all $\inn_{\tau} (f)$ for $f \in V$.
The following result and its proof are analogous to \cite[Corollary 1.6]{CO03}.

\begin{lem}
\label{helperdim} Let $\tau$ and $\sigma$ be term orders on $E$ and let $V \subset E_i$
be any $d$-dimensional subspace. If $m_1 \wedge \ldots \wedge m_d$  and $n_1 \wedge
\ldots \wedge n_d$ are $\tau$-standard monomials such that $m_1,\ldots,m_d $ is a $K$-basis of $\gin_{\tau} (V) $ and
$n_1,\ldots,n_d$ is a $K$-basis of $\gin_{\tau} (\inn_{\sigma} (V)) $, then $m_i \geq_{\tau} n_i$ for all $i = 1,\ldots,d$.
\end{lem}

Now we get analogously to the symmetric case \cite[Theorem 5.1]{CO04}:

\begin{thm}
\label{mainginleq}
Let $J \subset E$ be a graded ideal and $\tau$ a term order on $E$.
Then
$$
\beta_{ijp}^E(E/\gin(J)) \leq \beta_{ijp}^E(E/\gin_\tau(J)) \quad
\text{ for all } i,j,p.
$$
\end{thm}
\begin{proof}
Set $J'=\gin(J)$ and  $J''=\gin_\tau(J)$. Note  that $J'$ and $J''$
are squarefree
strongly stable ideals with the same Hilbert function as $J$. Thus by
Proposition \ref{maintechnical},
it is enough to show that $m_{\leq i}(J''_j) \leq m_{\leq i}(J'_j)$
for all $i$ and $j$.

Let $a_1,\dots,a_k$ be the generators of $J'_j$ and let $b_1\dots,b_k$
be those of
$J''_j$. We may assume that  the $a_r$'s and the $b_r$'s  are listed in the revlex order.
Then we claim:
$$
a_r\geq b_r
\text{ in the revlex order for all } r.
$$
This implies that $\max(a_r)\leq \max(b_r)$ for all $r$, and hence
$m_{\leq i}(J''_j)\leq
m_{\leq i}(J'_j)$. Thus it remains to shows the claim.

Let $V,W \subset E_j$ be two monomial vector spaces of the same dimension where $V$ has a
monomial $K$-basis $v_1,\dots,v_k$ and $W$ has a monomial $K$-basis $w_1,\dots,w_k$.
Define $V \geq_{\revlex} W$ if $v_i \geq_{\revlex} w_i$, $v_i>_{\revlex}
v_{i+1}$ and $w_i>_{\revlex} w_{i+1}$ for all $i$. Then Lemma \ref{helperdim} provides:
$$
\gin_{\revlex}(V)
\geq_{\revlex}
\gin_{\revlex}(\ini_\tau(V)).
$$
If $V \subset E_j$ is a generic subspace, then  $\ini_\tau(V)= \gin_\tau(V)$. Since
$\gin_{\revlex}(\gin_\tau(V))=\gin_\tau(V)$, we then get  $ \gin_{\revlex}(V)
\geq_{\revlex} \gin_\tau(V). $ Choosing $V=J_j'$ proves the claim.
\end{proof}

Observe that inequalities similar to the crucial
$m_{\leq i}((\gin_\tau(J))_j) \leq m_{\leq i}(\gin(J)_j)$ above
were shown by \cite[Proposition 2.4]{M05}.

\section{Simplicial complexes and algebraic shifting}
\label{sec-def}

In the second part of this paper we give some combinatorial applications
of the exterior algebra methods presented in Section \ref{cartan}
and some symmetric algebra methods which will be presented below.
For this we introduce at first some notation and discuss some basic concepts.
Recall that $\Gamma$ is called a {\em simplicial complex} on the
vertex set $[n]=\{1,\dots,n\}$ if $\Gamma$ is a subset of the power set of $[n]$ which is
closed under inclusion, i.~e.\ if $F \subseteq G$ and $G \in \Gamma$, then $F\in \Gamma$.
The elements $F$ of $\Gamma$ are called {\em faces}, and the maximal elements under
inclusion are called {\em facets}. We denote the set of all facets by $\facets(\Gamma)$.
If $F$ consists of $d+1$ elements of $[n]$, then $F$ is called a {\em $d$-dimensional} face,
and we write $\dim F = d$. The empty set is a face of dimension $-1$. Faces of dimension
0, 1 are called {\em vertices} and {\em edges}, respectively. $\Gamma$
is called {\em
pure} if all facets have the same dimension. If $\Gamma$ is non-empty,
then the {\em
dimension} $\dim \Gamma $ is the maximum of the dimensions of the
faces of $\Gamma$. Let
$f_i$ be the total number of $i$-dimensional faces of $\Gamma$ for
$i=-1,\dots,\dim
\Gamma$. The vector $f(\Gamma)=(f_{-1},\dots,f_{\dim \Gamma})$ is
called the {\em
$f$-vector} of $\Gamma$.

Several constructions associate other simplicial complexes to a given
one. For example,
the {\em Alexander dual} of $\Gamma$ is defined as $\Gamma^* =\{F
\subseteq [n] : F^c
\not\in \Gamma\}$ where $F^c=[n] \setminus F$. This gives a duality
operation on
simplicial complexes in the sense that $\Gamma^*$ is indeed a
simplicial complex on the
vertex set $[n]$ and we have that $(\Gamma^*)^*=\Gamma$.

The connection to algebra is due to the following construction. Let $K$ be a field and
$S=K[x_1,\dots,x_n]$ be the polynomial ring in $n$ indeterminates. $S$ is a graded ring
by setting $\deg x_i=1$ for $i=1,\dots,n$. We define the {\em Stanley--Reisner ideal}
$I_\Gamma=(\prod_{i \in F} x_i : F\subseteq [n],\ F\not\in \Gamma)$ and the corresponding
{\em Stanley--Reisner ring} $K[\Gamma]=S/I_\Gamma$. For a subset $F \subset [n]$, we
write $x_F$ for the squarefree monomial $\prod_{i \in F} x_i$. We will say that $\Gamma$
has an algebraic property like Cohen--Macaulayness if and only if $K[\Gamma]$ has this
property. Note that these definitions may depend on the ground field $K$.

There is an analogous exterior algebra construction due to Aramova, Herzog and Hibi (e.g.\ see \cite{HE01}).
Let $E=K \langle e_1,\dots,e_n\rangle$ be the exterior algebra on $n$ exterior variables.
$E$ is also a graded ring by setting $\deg e_i=1$ for $i=1,\dots,n$.
One defines the {\em exterior Stanley--Reisner ideal}
as $J_\Gamma=(e_F : F \not\in \Gamma)$
and the {\em exterior face ring} $K\{\Gamma\}=E/J_\Gamma$.

We keep these notations throughout the paper. Further information on simplicial complexes
can be found, for example,  in the books \cite{BRHE98} and
\cite{ST96}.

Let $\mathcal{C}([n])$ be the set of simplicial complexes on $[n]$. Following
constructions of Kalai, we define axiomatically the concept of algebraic shifting. See
\cite{HE01} and \cite{KA01} for surveys on this subject. We call a map $\shift \colon
\mathcal{C}([n]) \to \mathcal{C}([n])$ a {\em shifting operation} if the
following conditions are satisfied:

\begin{enumerate}
\item[(S1)]
If $\Gamma \in \mathcal{C}([n])$,
then $\shift(\Gamma)$ is a {\em shifted complex},
i.~e.\ for all $F \in \shift(\Gamma)$ and $j > i \in F$
we have that $(F  \setminus \{i\}) \cup \{j\} \in \shift(\Gamma)$.
\item[(S2)]
If $\Gamma \in \mathcal{C}([n])$ is shifted,
then $\shift(\Gamma)=\Gamma$.
\item[(S3)]
If $\Gamma \in \mathcal{C}([n])$, then
$f(\Gamma)=f(\shift(\Gamma))$.
\item[(S4)]
If $\Gamma', \Gamma \in \mathcal{C}([n])$ and $\Gamma' \subseteq \Gamma$ is a subcomplex,
then $\shift(\Gamma')$ is a subcomplex of $ \shift(\Gamma)$.
\end{enumerate}

See \cite{BNT2} for variations of these axioms.
Note that a simplicial complex $\Gamma$ is shifted
if and only if the Stanley--Reisner ideal
$I_\Gamma \subset S=K[x_1,\dots,x_n]$
is a {\em squarefree strongly stable ideal} with respect to $x_1>\dots>x_n$,
i.~e.\ for all $F\subseteq [n]$ with $x_F \in I_\Gamma$
and all $i \in F$, $j<i$, $j \not\in F$
we have that $(x_jx_F)/x_i \in I_\Gamma$.
Analogously,
$\Gamma$ is shifted
if and only if the exterior Stanley--Reisner ideal $J_\Gamma \subset
E=K\langle e_1,\dots,e_n \rangle$
is a {\em squarefree strongly stable ideal} with respect to $e_1>\dots>e_n$.

Two of the most important examples of such operations are defined as
follows.

\begin{exa}{(Symmetric algebraic shifting)}
\label{exsymmetric} At first we introduce the symmetric algebraic shifting
introduced in \cite{KA84}, \cite{KA91}. Here, we follow the algebraic approach of Aramova,
Herzog and Hibi \cite{ARHEHI00}.

Assume that $K$ is a field of characteristic $0$ and let $S=K[x_1,\dots,x_n]$.
We consider the following operation on monomial ideals of $S$.
For a monomial
$m=x_{i_1}\cdots x_{i_t}$ with $i_1\leq i_2 \leq \dots \leq i_t$
of $S$ we set $m^\sigma=x_{i_1}x_{i_2+1}\dots x_{i_t+ t-1}$.
For a monomial ideal $I$ with unique minimal system of generators $G(I)=\{m_1,\dots,m_s\}$
we set $I^\sigma=(m_1^\sigma,\dots,m_s^\sigma)$ in a suitable polynomial ring with sufficiently many variables.

Let $\Gamma$ be a simplicial complex on the vertex set $[n]$ with
Stanley--Reisner
ideal
$I_{\Gamma}$.
The {\em symmetric algebraic shifted complex} of $\Gamma$ is the unique
simplicial complex $\Delta^s(\Gamma)$ on the vertex set $[n]$ such that
$$
I_{\Delta^s(\Gamma)} = \bigl(\gin(I_\Gamma)\bigr)^\sigma \subset S.
$$

It is not obvious that $\Delta^s(\cdot)$ is indeed a shifting operation. The first
difficulty is to show that $I_{\Delta^s(\Gamma)}$ is an ideal of $S$. This and
the proofs of the other properties can be found in \cite{ARHEHI00} or \cite{HE01}.
\end{exa}

\begin{exa}{(Exterior algebraic shifting)}
\label{exexterior} Exterior algebraic shifting was also defined by Kalai in \cite{KA84}.
 Let $E=K\langle e_1,\dots,e_n\rangle$ where $K$ is any infinite field. Let
$\Gamma$ be a simplicial complex on the vertex set $[n]$. The {\em exterior algebraic
shifted complex} of $\Gamma$ is the unique simplicial complex $\Delta^e(\Gamma)$ on the
vertex set $[n]$ such that
$$
J_{\Delta^e(\Gamma)}
=
\gin(J_{\Gamma}).
$$
For an introduction to the theory
of Gr\"obner bases in the exterior algebra see \cite{ARHEHI97}. As opposed to symmetric algebraic shifting, it is much easier to see that $\Delta^ e(\cdot)$ is indeed a shifting operation
since generic initial ideals in the exterior algebra are already squarefree strongly
stable. See again \cite{HE01} for details.
\end{exa}

There are several other shifting operations. Since the proof of \cite[Proposition
8.8]{HE01} works for arbitrary term orders $\tau$, one can take  a generic initial ideal
$\gin_\tau(\cdot)$ in $E$ with respect to any term order $\tau$ to obtain, analogously to
Example \ref{exexterior}, the {\em exterior algebraic shifting operation} $\Delta^\tau(\cdot)$ with respect to the
term order $\tau$.


\section{Degree functions of simplicial complexes I}
\label{sec-degree}

In commutative algebra degree functions are designed to provide measures for the size and
the complexity of the structure of a given module. Given a simplicial complex $\Gamma$ on
the vertex set $[n]$,  we recall the definition of several degree functions on the
Stanley--Reisner ring $K[\Gamma]$ in terms of combinatorial data. The first important
invariant is the degree (or multiplicity) $\deg \Gamma$ of $\Gamma$. By definition $\deg
\Gamma=\deg K[\Gamma]$ is the degree of the Stanley--Reisner ring of $\Gamma$. It can be
combinatorially described as
\begin{equation} \label{eq-comb-deg}
\deg \Gamma=f_{\dim \Gamma}(\Gamma), \; \text{ the number of faces of maximal dimension.}
\end{equation}
Since any algebraic shifting operation $\shift(\cdot)$
preserves the $f$-vector,
we have that
\begin{equation} \label{eq-deg-inv}
\deg \Gamma = \deg \shift(\Gamma).
\end{equation}
In case $\Gamma$ is a Cohen--Macaulay complex, this invariant gives a lot of information
about $\Gamma$, because $\Gamma$ is pure and thus $\deg \Gamma$ counts all facets. There
are several other degree functions which take into account more information about
$\Gamma$.

Next we study the {\em arithmetic degree}.
Algebraically it is defined as
$$
\adeg \Gamma =
\adeg K[\Gamma]
=
\sum l(H^0_{\p}(K[\Gamma]_{\p})) \cdot  \deg \bigl(K[\Gamma]/\p \bigr)
$$
where the sum runs over the associated prime ideals of $K[\Gamma]$
and $l(\cdot)$ denotes the length function of a module.
In \cite{VA-98} and \cite{BNT} it was noted that the results of
\cite{STTRVO} imply that
\begin{equation} \label{eq-comb-adeg}
\adeg \Gamma = |\{F \in \facets(\Gamma)\}|
\end{equation}
is the number of facets of $\Gamma$.
Similarly, one defines for $i=0,\dots,\dim K[\Gamma]$
$$
\adeg_i \Gamma =
\adeg_i K[\Gamma]
=
\sum l(H^0_{\p}(K[\Gamma]_{\p})) \cdot  \deg \bigl(K[\Gamma]/\p \bigr)
$$
where the sum runs over the associated prime ideals of $K[\Gamma]$
such that $\dim K[\Gamma]/\p=i$ and gets
\begin{equation} \label{eq-comb-adeg-i}
\adeg_i \Gamma = |\{F \in \facets(\Gamma): \dim F=i-1\}|.
\end{equation}

By definition, shifting preserves the $f$-vector and the number of
facets of maximal
dimension of a simplicial complex. However, the total number of facets
may change under
shifting. But it may only increase, as we show now.

\begin{thm}
\label{thm-adeg-incr}
Let $\Delta$ be a shifting operation and
$\Gamma$ a $(d-1)$-dimensional simplicial complex. Then
$$
\adeg_i \Gamma \leq \adeg_i \Delta (\Gamma)
\text{ for } i=0,\dots,d.
$$
In particular,
$|\facets (\Gamma)|\leq |\facets (\Delta (\Gamma))|$.
\end{thm}

\begin{proof}
Let $m-1$ be the smallest dimension of a facet of
$\Gamma$. We prove the assertion by induction on the number $(d-1)-(m-1)=d-m \geq 0$.
If $d-m = 0$, then $\Gamma$ is a pure complex and only
$\adeg_d \Gamma$ is different from zero.
Using Identity \eqref{eq-deg-inv} we get
$\adeg_d \Gamma = \deg \Gamma = \deg \shift (\Gamma) = \adeg_d \shift (\Gamma)$ as claimed.

Let $d > m$. Observe that $\adeg_i \Gamma=0$ for $i<m$.
Denote by $F_1,\ldots,F_s \in \Gamma$ the $m-1$-dimensional
facets of
$\Gamma$ and let $\Gamma'$ be the subcomplex of $\Gamma$ whose facets are the remaining
facets of $\Gamma$. Thus, we have in particular:
$$
f_i (\Gamma') = f_i (\Gamma) \quad \fall i > m-1.
$$
Since shifting does not change the $f$-vector and $\shift (\Gamma')$ is a subcomplex of
$\shift (\Gamma)$, we conclude that
$$
f_i (\shift (\Gamma')) = f_i (\shift (\Gamma)) \quad \fall i > m-1,
$$
and, in particular, that each facet of $\shift (\Gamma')$ of dimension $> m-1$ is also a
facet of $\shift (\Gamma)$.
It follows from the induction hypothesis that
$$
\adeg_i \Gamma
=
\adeg_i \Gamma'
\leq
\adeg_i \Delta(\Gamma')
=
\adeg_i \Delta(\Gamma)
\text{ for }
i>m.
$$
Hence, it remains to show that
$\shift (\Gamma)$ has at least $s=\adeg_{m} \Gamma$ facets of dimension $m-1$.

To this end note that, by definition of $\Gamma'$, $f_{m-1} (\Gamma) = f_{m-1} (\Gamma') + s$,
thus we get for the shifted complexes $f_{m-1} (\shift (\Gamma)) = f_{m-1} (\shift (\Gamma')) +
s$. Let $G \in \shift (\Gamma)$ be any $(m-1)$-dimensional face that is not in
$\shift (\Gamma')$. Assume that $G$ is strictly contained in a face $\widetilde{G}$ of
$\shift (\Gamma)$. Then $\dim \widetilde{G} > m-1$ and, using the fact that the faces of
$\shift (\Gamma')$ and $\shift (\Gamma)$ of dimension $> m-1$ coincide, we conclude that
$\widetilde{G} \in \shift (\Gamma')$.
Thus $G \in \shift (\Gamma')$. This contradiction
to the choice of $G$ shows that $G$ is a facet of $\shift (\Gamma)$ and we are done.
\end{proof}

Observe that we only used axioms (S3) and (S4) of a shifting operation to prove
Theorem \ref{thm-adeg-incr}.
It would be interesting to know if there are further results in this direction like:
\begin{quest}
\label{importantquestion}
Is there a shifting operation $\shift(\cdot)$ such that $\adeg$ increases the least,
i.~e.\, for every simplicial complex $\Gamma$, we have $\adeg \shift (\Gamma) \leq \adeg
\shift' (\Gamma)$ for each algebraic shifting operation $\shift'(\cdot)$?
\end{quest}

We use the results in Section \ref{cartan} to give a partial answer to the previous question.
We begin by observing that shifting and Alexander duality commute.

\begin{lem}
\label{shifthelper} Let $\Gamma$ be a simplicial complex on the vertex
set $[n]$. Let
$\Delta^{\tau}(\cdot)$ be the  exterior shifting operation with respect to a
term order $\tau$ on $E$. Then Alexander duality commutes with shifting, i.\ e.\ we have $
\Delta^\tau(\Gamma)^* = \Delta^\tau(\Gamma^*). $
\end{lem}
\begin{proof}
In \cite{HETE} this was proved for $\tau$ being the revlex order. But the proof works for
any term order $\tau$.
\end{proof}

As a further preparation, we need an interpretation of the arithmetic degree over the
exterior face ring using the socle. The socle $\Soc N$ of a finitely generated $E$-module
$N$ is the set of elements $x \in N$ such that $(e_1,\dots,e_n)x=0$. Observe that $\Soc
N$ is always a finite-dimensional $K$-vector space.

\begin{lem}
\label{arithdesc} Let $\Gamma$ be a $(d-1)$-dimensional
simplicial complex on the vertex set $[n]$.
Then
$$
\adeg_i  \Gamma = \dim_K \Soc K\{\Gamma\}_i
\text{ for }
i=0,\dots,d.
$$
In particular, $\adeg  \Gamma = \dim_K \Soc K\{\Gamma\}$.
\end{lem}

\begin{proof}
The residue classes of the monomials $e_F$, $F \in \Gamma$, form a $K$-vector space basis
of $K\{\Gamma\}$. Hence, the facets of $\Gamma$ correspond to a $K$-basis of $\Soc
K\{\Gamma\}$.
\end{proof}

We have seen in Theorem \ref{thm-adeg-incr} that each shifting operation increases  the
number of facets. Below we show that among the exterior shifting operations, standard
exterior shifting with respect to revlex leads to the least possible increase.

Now we recall the concept of a sequentially Cohen--Macaulay module. Let $M$ be a finitely
generated graded $S$-module. The module $M$ is said to be {\em sequentially
Cohen--Macaulay} (sequentially CM modules for short), if there exists a finite filtration
\begin{equation}
\label{filter}
0=M_0 \subset M_1 \subset M_2 \subset \dots \subset M_r=M
\end{equation}
of $M$ by graded submodules of $M$
such that each quotient $M_i/M_{i-1}$ is Cohen--Macaulay and
$\dim M_1/M_0 <\dim M_2/M_1 <\dots<\dim M_r/M_{r-1}$ where $\dim$ denotes
the Krull dimension of $M$.

\begin{thm}
\label{shifted}
Let $\Gamma$ be a $(d-1)$-dimensional simplicial complex on the vertex set $[n]$. For any term order $\tau$ on $E$, we have
$$
\adeg_i \Gamma \leq \adeg_i \Delta^e(\Gamma) \leq \adeg_i \Delta^\tau(\Gamma)
\text{ for }
i=0,\dots,d.
$$
In particular,
$\adeg\Gamma \leq \adeg\Delta^e(\Gamma) \leq \adeg \Delta^\tau(\Gamma)$.
Moreover, $\adeg \Gamma = \adeg \Delta^e(\Gamma)$ if and only if $\Gamma$ is sequentially
Cohen--Macaulay.
\end{thm}

\begin{proof}
Using $\adeg_i \Gamma = \dim_K \Soc (E/J_{\Gamma})_i$ and the definition of Alexander duality,
it is easy to see that $\adeg_i \Gamma$ coincides with the number of minimal generators of
$J_{\Gamma^*}$ of degree $n-i$ which is $\beta_{0n-i}^E(J_{\Gamma^*})$,
i.~e.\ we have
$ \adeg_i \Gamma = \beta_{0n-i}^E(J_{\Gamma^*}) =
\beta_{1n-i}^E(E/J_{\Gamma^*}).$ Hence Lemma \ref{shifthelper} provides
$$
\adeg_i \Gamma
=
\beta_{1n-i}^E(E/J_{\Gamma^*})
\leq
\beta_{1n-i}^E(E/\gin(J_{\Gamma^*}))
=
\beta_{1n-i}^E(E/J_{\Delta^e(\Gamma^*)})
$$
$$
= \adeg_i \Delta^e(\Gamma^*)^* = \adeg_i \Delta^e(\Gamma)
$$
for $i=0,\dots,d$.
Since $ \beta_{1j}^E(E/J_{\Gamma^*}) \leq
\beta_{1j}^E(E/\gin(J_{\Gamma^*})) $ for all
$j$, it follows from Theorem \ref{maincl} that we have equality if and only if
$J_{\Gamma^*}$ is
componentwise linear. By \cite{HEHI99},  this is the case if and only if the
Stanley--Reisner ideal $I_{\Gamma^*} \subset S=K[x_1,\dots,x_n]$ is
componentwise linear.
But the latter is equivalent to $\Gamma$ being  sequentially
Cohen--Macaulay according to
\cite{HEREWE}.

Combining the above argument and  Theorem \ref{mainginleq}, we obtain also
$$
\adeg_i \Delta^e(\Gamma) = \beta_{1n-i}^E(E/\gin(J_{\Gamma^*})) \leq
\beta_{1n-i}^E(E/\gin_\tau(J_{\Gamma^*})) = \adeg_i \Delta^\tau(\Gamma),
$$
and the proof is complete.
\end{proof}

\begin{rem}
\
\begin{enumerate}
\item
Observe that Theorem \ref{mainginleq} is essentially equivalent to the
following inequalities $m_{\leq i}((\gin_\tau(J))_j) \leq m_{\leq i}(\gin(J)_j)$.
Thus
one could argue in the above proof of Theorem \ref{shifted}
by using the latter numbers instead of the corresponding graded Betti-numbers.
(Again see also \cite[Proposition 2.4]{M05} for the needed inequalities.)
\item
The fact that
$\adeg \Gamma = \adeg \Delta^e(\Gamma)$ if and only if $\Gamma$ is sequentially
Cohen--Ma\-caulay can also be shown by refining the arguments in the proof of Theorem \ref{thm-adeg-incr}
and using the following facts:
a simplicial complex $\Gamma$ is sequentially Cohen--Macaulay
if and only if each subcomplex $\Gamma^{[i]}$ generated by the $i$-faces of $\Gamma$
is Cohen--Macaulay (see \cite{DU});
$\Gamma$ is Cohen--Macaulay if and only if $\Delta^e(\Gamma)$ is pure.
\end{enumerate}
\end{rem}

The arithmetic degree has nice properties especially for Cohen--Macaulay $K$-algebras.
However, there are some disadvantages for non-Cohen--Macaulay $K$-algebras. See
\cite{VA} for a discussion. Vasconcelos axiomatically
defined the following
concept. Recall that $S=K[x_1,\dots,x_n]$ where $K$ is a field. A numerical function
$\Deg$ that assigns to every finitely generated graded $S$-module  a non-negative integer
is said to be an {\em extended degree function} if it satisfies the following conditions:
\begin{enumerate}
\item
If $L=H^0_{\m}(M)$, then $\Deg M =\Deg M/L +l(L).$
\item
If $y\in S_1$ is sufficiently general and $M$-regular,
then $\Deg M \geq \Deg M/yM.$
\item
If $M$ is a Cohen--Macaulay module, then $\Deg M=\deg M.$
\end{enumerate}
There exists a {\em smallest extended degree function} $\sdeg$ in the sense that $\sdeg M
\leq \Deg M$ for every other extended degree function.

\begin{rem}
\label{sdeghelper}
We need the following properties of $\sdeg$. (For proofs see
\cite{NR}.) For simplicity we state them only for a graded $K$-algebra $S/I$ where $I
\subset S$ is a graded ideal.
\begin{enumerate}
\item
$S/I$ is Cohen--Macaulay if and only if $\sdeg S/I= \deg S/I$. (This is true for any
extended Degree function.)
\item
$\adeg S/I\leq \sdeg S/I$.
\item
If $S/I$ is sequentially Cohen--Macaulay, then $\adeg S/I=\sdeg S/I$.
\item
$\sdeg S/I = \sdeg S/\gin(I)$.
\end{enumerate}
\end{rem}

Then we have:

\begin{thm}
\label{sedgadegshift}
Let $\Gamma$
be a simplicial complex on the vertex set $[n]$.
We have that:
\begin{enumerate}
\item
$\sdeg \shift(\Gamma)=\adeg \shift(\Gamma)$ for each algebraic shifting operation
$\shift$.
\item
$\deg \Gamma \leq \adeg \Gamma \leq
\sdeg \Gamma = \sdeg \shift^s(\Gamma)=\adeg \shift^s(\Gamma)$.
\end{enumerate}
\end{thm}
\begin{proof}
In the proof we use that $\shift(\Gamma)$ is always sequentially Cohen--Macaulay, i.\ e.\
$K[\shift(\Gamma)]$ is a sequentially CM-ring (cf.\ Proposition \ref{prop-shifted-sCM}
below).  This fact and the properties of $\sdeg$ observed in Remark \ref{sdeghelper}
imply (i).

The only critical equality in (ii) is
$\sdeg \Gamma = \sdeg \shift^s(\Gamma)$.
We compute
\begin{eqnarray*}
\sdeg \Gamma
&=& \sdeg S/I_\Gamma\\
&=& \sdeg S/\gin(I_\Gamma)\\
&=& \adeg S/\gin(I_\Gamma)\\
&=& \adeg S/I_{\shift^s(\Gamma)}\\
&=& \sdeg S/I_{\shift^s(\Gamma)}
\end{eqnarray*}
where the first equality is the definition of
$\sdeg \Gamma$,
the second one was noted above and
the third equality follows from the fact that $S/\gin(I_\Gamma)$
is sequentially Cohen--Macaulay (see \cite[Theorem 2.2]{HESB}).
The forth equality
follows again from
\cite[Theorem 6.6]{BNT}
and the last equality is (i) which we proved already.
\end{proof}

We already mentioned
that $\deg \Gamma$ is the number of facets of maximal dimension of $\Gamma$
and
$\adeg \Gamma$ is the number
of facets of $\Gamma$.
Theorem \ref{sedgadegshift}
provides a combinatorial interpretations of $\sdeg \Gamma$.
It is the number of facets of $\Delta^s(\Gamma)$.


\section{The Cohen--Macaulay property and iterated Betti numbers}
\label{ringprop}

In this section we use degree functions to relate properties of the shifted complex and
iterated Betti numbers to the original complex.

Our first observation is well-known to specialists. It states that shifted complexes have
a nice algebraic structure.

\begin{prop}
\label{prop-shifted-sCM}
Shifted simplicial complexes are sequentially Cohen--Macaulay.

In particular,
if $\Gamma$ is a simplicial complex on the vertex set $[n]$
and $\shift(\cdot)$ is an arbitrary algebraic shifting operation,
then $\shift(\Gamma)$ is sequentially Cohen--Macaulay.
\end{prop}

\begin{proof}
Let $\Gamma$ be a shifted simplicial complex.
Recall that then $I_{\Gamma} \subset S$ is a squarefree strongly stable ideal
with respect to $x_1>\dots>x_n$, i.~e.\ for all
$x_F=\prod_{l \in F}x_l \in I_{\Gamma}$
and $i$ with $x_i|x_F$ we have for all $j<i$ with $x_j\nmid x_F$
that $(x_F/x_i)x_j\in I_{\Gamma}$.

It is easy to see that we have that $I_{\Gamma}$ is squarefree strongly stable if and
only if $I_{\Gamma^*}$ is squarefree strongly stable where $\Gamma^*$ is the Alexander
dual of $\Gamma$. It is well-known  that in this situation $I_{\Gamma^*}$ is a so-called
componentwise linear ideal (see Section \ref{cartan} for the definition) and thus the
theorem follows now from \cite[Theorem 9]{HEREWE}.
\end{proof}

\begin{rem}
Alternatively one
can prove Proposition \ref{prop-shifted-sCM}
using more combinatorial arguments as follows.
Shifted simplicial complexes are non-pure shellable by \cite{BW97}
and thus sequentially Cohen--Macaulay by \cite{ST96}.
\end{rem}

It is easy to decide whether a shifted simplicial complex is Cohen--Macaulay because of
the following well-known fact:

\begin{prop}
Let $\Gamma$ be a shifted simplicial complex on the vertex set $[n]$. Then the  following
statements are equivalent:
\begin{enumerate}
\item
$\Gamma$ is Cohen--Macaulay;
\item
$\Gamma$ is pure.
\end{enumerate}
\end{prop}

\begin{proof}
If $\Gamma$ is Cohen--Macaulay, then it is well-known that
$\Gamma$ is pure.
Assume now that $\Gamma$ is a pure shifted simplicial complex.
We compute
$$
\sdeg S/I_{\Gamma}
=
\adeg S/I_{\shift^s(\Gamma)}
=
\adeg S/I_{\Gamma}
=
\deg S/I_{\Gamma}.
$$
The first equality was shown in Theorem \ref{sedgadegshift}, the second one holds because
$\Gamma$ is shifted, and the third one because  $\Gamma$ is pure. Hence $\sdeg
S/I_{\Gamma}=\deg S/I_{\Gamma}$ and this implies by Remark \ref{sdeghelper} that $\Gamma$
is Cohen--Macaulay.
\end{proof}

Intuitively, shifting leads to a somewhat simpler complex and one would like to transfer
properties from $\shift (\Gamma)$ to $\Gamma$. In this respect,  we have:

\begin{prop}
\label{justnoted} Let  $\Gamma$ be a simplicial complex on the vertex set $[n]$ and let
$\shift(\cdot)$ be any algebraic shifting operation such that $\sdeg \Gamma \leq \sdeg
\shift(\Gamma)$. If $\shift(\Gamma)$ is Cohen--Macaulay, then $\Gamma$ is
Cohen--Macaulay.
\end{prop}
\begin{proof}

Suppose that $\shift(\Gamma)$ is Cohen--Macaulay. Then $\shift(\Gamma)$ is pure and we
have that $\adeg \shift(\Gamma) = \deg \shift(\Gamma)$. Thus, we get
$$
\deg \Gamma
\leq
\sdeg \Gamma
\leq
\sdeg \shift(\Gamma)
=
\adeg \shift(\Gamma)
=
\deg \shift(\Gamma)
=
\deg \Gamma.
$$
The only critical equality is
$\sdeg \shift(\Gamma)
= \adeg \shift(\Gamma)$ which follows from Remark \ref{sdeghelper} (iii) and Proposition \ref{prop-shifted-sCM}.
Hence $\deg \Gamma=\sdeg \Gamma$ and it follows that $\Gamma$ is Cohen--Macaulay.
\end{proof}

For exterior algebraic shifting,  the following result was proved in \cite{KA01}.
For symmetric algebraic shifting, this result follows from some general homological
arguments given in \cite{BYCHPO}, as noted in \cite{BNT} and essentially first appeared in \cite{KA91} (Theorem 6.4). Below, we provide a  very short proof using only degree functions.
\begin{thm}
\label{cmnice} Let  $\Gamma$ be a simplicial complex on the vertex set $[n]$. Then the
following statements are equivalent:
\begin{enumerate}
\item
$\Gamma$ is Cohen--Macaulay;
\item
$\shift^{s}(\Gamma)$ is Cohen--Macaulay;
\item
$\shift^{s}(\Gamma)$ is pure.
\end{enumerate}
\end{thm}
\begin{proof}
We already know that (ii) and (iii) are equivalent.
Suppose that (i) holds.
Then
$$
\deg \shift^s(\Gamma)
=
\deg \Gamma
=
\sdeg \Gamma
=
\adeg \shift^s(\Gamma)
\geq
\adeg \Gamma
\geq
\deg \Gamma
=
\deg \shift^s(\Gamma).
$$
The first equality is trivial.
The second one holds because $\Gamma$ is Cohen--Macaulay,
and the third one follows from Theorem \ref{sedgadegshift} (ii).
The first inequality is due to
Theorem \ref{thm-adeg-incr} and the second one follows from the definition of $\adeg$.
Hence $\adeg \shift^s(\Gamma) = \deg \shift^s(\Gamma)$, and thus $\shift^s(\Gamma)$ is
pure, as claimed in (iii).

The fact that (iii) implies (i) follows from   $\sdeg \Gamma=\adeg \shift^s(\Gamma) \leq
\sdeg \shift^s (\Gamma)$ and Proposition \ref{justnoted}.
\end{proof}

\begin{rem}
Using Theorem \ref{shifted} one can prove in a similar way that Theorem \ref{cmnice} holds also for the exterior shifting $\Delta^e$.
\end{rem}

For the next results, we need some further notation and definitions. At first we recall a
refinement of the $f$- and $h$-vector of a simplicial complex due to \cite{BW96}.

\begin{defi}
Let $\Gamma$ be a $(d-1)$-dimensional simplicial complex on the vertex set $[n]$. We
define
$$
f_{i,r}(\Gamma) =  |\{ F \in \Gamma : \deg_{\Gamma} F = i \text{ and } |F| = r \}|
$$
where $\deg_{\Gamma} F= \max \{|G| : F \subseteq G \in \Gamma\}$ is the {\em degree} of a
face $F \in \Gamma$. The triangular integer array $f(\Gamma) = (f_{i,r}(\Gamma))_{0\leq r
\leq i\leq d}$ is called the {\em $f$-triangle} of $\Gamma$. Further we define
$$
h_{i,r}(\Gamma) =
\sum_{s=0}^r
(-1)^{r-s}\binom{i-s}{r-s} f_{i,s}(\Gamma).
$$
The triangular array $h(\Gamma)= (h_{i,r}(\Gamma))_{0\leq r\leq i\leq d}$
is called the {\em $h$-triangle} of $\Gamma$.
\end{defi}
It is easy to see that giving the $h$-triangle, one can also compute the $f$-triangle and
thus these triangles determine each other. For $F \subseteq [n]$ let $\init(F)=\{k, \dots, n\}$ if $k,\dots,n \in F$, but $k-1 \not\in F$. (Here we set $\init(F)=\emptyset$ if no such $k$ exists.) Duval proved in \cite[Corollary 6.2]{DU} that if $\Gamma$ is sequentially
Cohen--Macaulay, then $h_{i,r} (\Gamma)= |\{F \in \facets(\Gamma) : |\init(F)|=i-r,\ |F|=i
\}|$. In particular, this formula holds for shifted simplicial complexes. Using this fact we define:

\begin{defi}
Let  $\Gamma$ be a simplicial complex on the vertex set $[n]$
and $\shift(\cdot)$ an arbitrary algebraic shifting operation.
We call
$$
b^\Delta_{i,r}(\Gamma)
=
h_{i,r}(\shift(\Gamma))
=
|\{F \in \facets(\shift(\Gamma)) : |\init(F)|=i-r,\ |F|=i \}|
$$
the {\em $\Delta$-iterated Betti numbers} of $\Gamma$.
\end{defi}
For symmetric and exterior shifting these
notions were defined and studied in \cite{BNT} and \cite{DURO}.
Since $\shift(\shift(\Gamma))=\shift(\Gamma)$
we have that
$$
b^\Delta_{i,r}(\Gamma)
=
b^\Delta_{i,r}(\shift(\Gamma))
$$
and therefore $\Delta$-iterated Betti numbers are invariant under each algebraic shifting
operation $\shift(\cdot)$. As noted in \cite{BNT},  there is the following comparison
result.
\begin{thm} \label{thm-it-Betti}
Let $\Gamma$
be a simplicial complex on the vertex set $[n]$
and $\Delta(\cdot)$ be the symmetric shifting $\Delta^s$ or the exterior shifting $\Delta^e $.
Then the following conditions are equivalent:
\begin{enumerate}
\item
$\Gamma$ is sequentially Cohen--Macaulay;
\item
$b^{\shift}_{i,r}(\Gamma)=h_{i,r}(\Gamma)$
for all $i,r$.
\end{enumerate}
In particular, if $\Gamma$ is sequentially Cohen--Macaulay,
then
all iterated Betti numbers
of $\Gamma$ with respect to the symmetric  shifting and the exterior
algebraic shifting coincide.
\end{thm}
\begin{proof}
Duval proved in \cite[Theorem 5.1]{DU} the equivalence of (i) and (ii) for exterior
algebraic shifting $\shift^e$. But his proof works also for the symmetric shifting 
since he used only axioms (S1)--(S4) of an algebraic shifting operation and the fact that Theorem \ref{cmnice} holds also for the exterior shifting.  

\end{proof}

The next  result  extends Theorem 6.6 of \cite{BNT} where the case of the symmetric
algebraic shifting operation was studied. It follows by combining Identities
\eqref{eq-comb-deg}, \eqref{eq-comb-adeg}, Proposition \ref{prop-shifted-sCM}, and Remark
\ref{sdeghelper} (iii).

\begin{thm}
Let $\Gamma$ be a simplicial complex on the vertex set $[n]$ of dimension $d-1$ and let
$\shift(\cdot)$ be any algebraic shifting operation. Then we have:
\begin{enumerate}
\item
$
\deg \shift(\Gamma) = \sum_{r=0}^{d} b^{\shift}_{d,r}(\Gamma)
= |\{F \in \facets(\shift(\Gamma)),\ \dim F=d-1 \}|.
$
\item
$
\sdeg \shift(\Gamma)= \adeg \shift(\Gamma)
= \sum_{r,i} b^{\shift}_{i,r}(\Gamma)
= |\{F \in \facets(\shift(\Gamma))\}|.
$
\end{enumerate}

\end{thm}

The last two results suggest the following question:

\begin{quest}
\ What is the relationship between $b^{\shift}_{i,r}(\Gamma)$ for different shifting
operations? For example,  in \cite{BNT} it is conjectured that
$b^{\shift^s}_{i,r}(\Gamma) \leq b^{\shift^e}_{i,r}(\Gamma)$ for all $i,r$.
\end{quest}

\section*{Acknowledgments}
The authors would like to thank the referee for the careful reading and the very helpful comments.




\end{document}